\documentclass[11pt]{amsart}
\usepackage{amscd,amssymb,verbatim}

\theoremstyle{plain}
\newtheorem{thm}{Theorem}

\newtheorem{lem}[thm]{Lemma}
\newtheorem{prop}[thm]{Proposition}

\newcommand{\ep}{\varepsilon}

\newcommand{\vp}{\varphi}

\newcommand{\bs}{\backslash}

\newcommand{\co}{\operatorname{co}}

\newcommand{\N}{{\Bbb N}} 
\newcommand{\R}{{\Bbb R}}

\newcommand{\supp}{\operatorname{supp}}
\newcommand{\sump}{\sum_{\vp\in A_i}}
\newcommand{\ki}{{^k_{i=1}}}

\begin{document}

\title[Schachermayer's space under renorming]
{An asymptotic property of Schachermayer's space
under renorming}
\author{Denka Kutzarova}
\address{Institute of Mathematics\\Bulgarian Academy of Sciences\\
Sofia\\ Bulgaria}
\curraddr{Department of Mathematics \\
         University of South Carolina \\
         SC 29208\\
         U.S.A.}
\email{denka@math.sc.edu}
\author{Denny H. Leung}
\address{Department of Mathematics \\
         National University of Singapore \\
         Singapore 117543}
\email{matlhh@nus.edu.sg}


\begin{abstract} 
Let $X$ be a Banach space with closed unit ball $B$. Given $k \in \N$, $X$ is said to be $k$-$\beta$, repectively, $(k+1)$-nearly uniformly convex ($(k+1)$-NUC), if for every $\ep > 0$, there exists $\delta$, $0 < \delta < 1$, so that for every $x \in B$, and every $\ep$-separated sequence $(x_n) \subseteq B$, there are indices $(n_i)^k_{i=1}$, respectively, $(n_i)^{k+1}_{i=1}$, such that
$ \frac{1}{k+1}\|x + \sum^{k}_{i=1}x_{n_i}\| \leq 1- \delta$,
respectively,
$\frac{1}{k+1}\|\sum^{k+1}_{i=1}x_{n_i}\| \leq 1- \delta$.
It is shown that a Banach space constructed by Schachermayer is $2$-$\beta$, but is not isomorphic to any $2$-NUC Banach space. Modifying this example, we also show that there is a $2$-NUC Banach space which cannot be equivalently renormed to be $1$-$\beta$.
\end{abstract}

\maketitle


In \cite{H}, R.\ Huff introduced the notion of nearly uniform convexity (NUC). A Banach space $X$ with closed unit ball $B$ is said to be NUC if for any $\ep > 0$, there exists $\delta < 1$ such that for every $\ep$-separated sequence in $B$, $\co((x_n)) \cap \delta B \neq \emptyset$. Here $\co(A)$ denotes the convex hull of a set $A$; a sequence $(x_n)$ is $\ep$-separated if $\inf\{\|x_n-x_m\| : m \neq n\} \geq \ep$. Huff showed that a Banach space is NUC if and only if it is reflexive and has the uniform Kadec-Klee property (UKK). Recall that a Banach space $X$ with closed unit ball $B$ is said to be UKK if for any $\ep > 0$, there exists $\delta < 1$ such that for every $\ep$-separated sequence $(x_n)$ in $B$ which converges weakly to some $x \in X$, we have $\|x\| \leq \delta$. A recent result of H.\ Knaust, E.\ Odell, and Th.\ Schlumprecht \cite{KOS} gives an isomorphic characterization of spaces having NUC. They showed that a separable reflexive Banach space $X$ is isomorphic to a UKK space if and only if $X$ has finite Szlenk index. 

Another property related to NUC is the property ($\beta$) introduced by Rolewicz \cite{R}. In \cite{K1}, building on the work of S.\ Prus \cite{P1,P2}, the first author showed that a separable Banach space $X$ is isomorphic to a space with ($\beta$) if and only if both $X$ and $X^*$ are isomorphic to NUC spaces. In \cite{K}, a sequence of properties lying in between ($\beta$) and NUC are defined. Let $X$ be a Banach space with closed unit ball $B$. Given $k \in \N$, $X$ is said to be $k$-$\beta$, repectively, $(k+1)$-nearly uniformly convex ($(k+1)$-NUC), if for every $\ep > 0$, there exists $\delta$, $0 < \delta < 1$, so that for every $x \in B$, and every $\ep$-separated sequence $(x_n) \subseteq B$, there are indices $(n_i)^k_{i=1}$, respectively, $(n_i)^{k+1}_{i=1}$, such that
\[ \frac{1}{k+1}\|x + \sum^{k}_{i=1}x_{n_i}\| \leq 1- \delta, \]
respectively,
\[ \frac{1}{k+1}\|\sum^{k+1}_{i=1}x_{n_i}\| \leq 1- \delta. \]
It follows readily from the definitions that every $k$-$\beta$ space is $(k+1)$-NUC, every $(k+1)$-NUC space is $(k+1)$-$\beta$, and that every $k$-$\beta$ space (or $(k+1)$-NUC space) is NUC. It is proved in \cite{K} that property $1$-$\beta$ is equivalent to the property ($\beta$) of Rolewicz. 
It is worth noting that the ``non-uniform'' version of property $k$-NUC has been well-studied. For $k \geq 2$, a Banach space $X$ is said to have property ($kR$) if every sequence $(x_n)$ in $X$ which satisfies $\lim_{n_1}\dots\lim_{n_k}\|x_{n_1}+\dots+x_{n_k}\| = k\lim_n\|x_n\|$ is convergent \cite{FG}. It is clear that property ($kR$) implies property ($(k+1)R$). It follows from James' characterization of reflexivity that every ($kR$) space is reflexive. A recent result of E.\ Odell and Th.\ Schlumprecht \cite{OS} shows that a separable Banach space is reflexive if and only if it can be equivalently renormed to have property ($2R$). Thus, all the properties ($kR$) are isomorphically equivalent.
Similarly, ``non-asymptotic'' properties known as $k$-uniform rotundity have been studied \cite{S}. These properites are also isomorphically equivalent to each other as they are all equivalent to superreflexivity. 
In this paper, we find that the situation is different for the properites $k$-NUC and $k$-$\beta$. To be precise, we use the space constructed by W.\ Schachermayer in \cite{Sch} and a variant to distinguish the properties $1$-$\beta$, $2$-NUC, and $2$-$\beta$ isomorphically.

Let $T = \cup^\infty_{n=0}\{0,1\}^n$ 
be the dyadic tree.  If $\vp = (\ep_i)^m_{i=1}$ and $\psi = (\delta_i)^n_{i=1}$ are nodes in $T$, we say that $\vp \leq \psi$ if $m \leq n$ and $\ep_i = \delta_i$ for $1 \leq i \leq m$. Also $\emptyset \leq \vp$ for all $\vp \in T$. Two nodes $\vp$ and $\psi$ are said to be comparable if either $\vp \leq \psi$ or $\psi \leq \vp$; they are incomparable otherwise.
Let $\vp \in T$, 
denote by $T_\vp$ or $T(\vp)$ the subtree rooted at 
$\vp$, i.e., the subtree consisting of all nodes $\psi$ 
such that $\vp \leq \psi$. A node $\vp \in T$ has length $n$ if $\vp \in \{0,1\}^n$. The length of $\vp$ is denoted by $|\vp|$. Given $\vp = (\ep_i)^n_{i=1} \in T$, let $S_\vp$ be the set consisting of all nodes $\psi = (\delta_i)^m_{i=1}$ such that $m \geq n$, $\delta_i = \ep_i$ if $1 \leq i \leq n$, and $\delta_i = 0$ otherwise.
Say that a subset $A$ of $T$ is admissible, respectively, acceptable, if there 
exists $n \in \N\ \cup\{0\}$ 
such that (a) $A \subseteq \cup_{|\vp|=n}T_\vp$, 
and (b) $|A \cap T_\vp| \leq 1$ for all 
$\vp$ with $|\vp| = n$, respectively, 
(a$'$) $A \subseteq \cup_{|\vp|=n}S_\vp$, 
and (b$'$) $|A \cap S_\vp| \leq 1$ for all 
$\vp$ with $|\vp| = n$. For subsets $A$ and $B$ of $T$, 
say that $A \ll B$ if 
$\max\{|\vp| : \vp \in A\} < \min\{|\vp| : \vp \in B\}$. 
Let $c_{00}(T)$ be the space
of all finitely supported real-valued functions defined on $T$. 
For $x \in c_{00}(T)$, 
let 
\[ \|x\|_X = \sup(\sum^k_{i=1}(\sum_{\vp\in A_i}|x(\vp)|)^2)^{1/2}, \]
where the sup is taken over all $k \in \N$, and all sequences 
of admissible subsets 
$A_1 \ll A_2 \ll \dots \ll A_k$ of $T$. The norm $\|\cdot\|_Y$ is defined similarly except that the sup is taken over all sequences of acceptable subsets $A_1 \ll A_2 \ll \dots \ll A_k$ of $T$.
Schachermayer's space $X$ 
is the completion of 
$c_{00}(T)$ with respect to the norm $\|\cdot\|_X$. The completion of $c_{00}(T)$ with respect to $\|\cdot\|_Y$ is denoted by $Y$.\\

\noindent{\bf Remark.} The space $X$ defined here differs 
from Schachermayer's original 
definition and is only isomorphic to the space defined in \cite{Sch}.\\

In \cite{K}, it was  
shown that $X$ (with the norm given in \cite{Sch}) is $8$-NUC but is not isomorphic to any $1$-$\beta$ space. We first show that $(X,\|\cdot\|_X)$ and $(Y, \|\cdot\|_Y)$ are $2$-$\beta$ and $2$-NUC respectively. We begin with a trivial lemma 
concerning the $\ell^2$-norm $\|\cdot\|_2$.

\begin{lem} \label{trivial}
If $\alpha$, $\beta$, and $\gamma$ are vectors in the unit
ball of $\ell^2$, and $\|\alpha + \beta + \gamma\|_2/3 \geq  1 - \delta$, 
then $\max\{\|\alpha - \beta\|_2, \|\alpha - \gamma\|_2, 
\|\beta - \gamma\|_2\} \leq \sqrt{18\delta}$.
\end{lem}

\begin{prop}
$(X,\|\cdot\|_X)$ is $2$-$\beta$.
\end{prop}

\begin{proof}
Let $x$, and $x_n, n \geq 1$ be elements in the unit 
ball of $X$ such that $(x_n)$ 
is $\ep$-separated. Choose $\delta > 0$ such that
\begin{equation}\label{delta}
(1 - 3\delta)^2 + [(1 - 24\delta)^{1/2} - (1 - \ep^2/9)^{1/2}]^2 > 1.
\end{equation}
%
%
%
Without loss of generality, we may assume that $(x_n)$ converges pointwise (as a sequence of functions on $T$) to some $y_0 : T \to \R$. It is clear that if $y, z \in X$ and $\supp y \ll \supp z$, then $\|y+z\|^2_X \geq \|y\|^2_X + \|z\|^2_X$. It follows easily that $y_0 \in X$. 
Let $y_n = x_n - y_0$.  It may be assumed that $(\|y_n\|_X)$ converges. As $(x_n)$ is $\ep$-separated, so is $(y_n)$. We may thus further assume that $\|y_n\|_X > \ep/3$ for all $n \in \N$. By going to a subsequence and perturbing
the vectors $x$, $y_0$ and 
$y_n, n \geq 1$ by as little as we please, it may be further assumed that 
(a) they all belong 
to $c_{00}(T)$,
(b) $\supp x \cup \supp y_0 \ll \supp y_1 \ll \supp y_2$, and 
(c) $\|y_1\chi_{T\vp}\|_\infty = \|y_2\chi_{T\vp}\|_\infty$ 
for all $\vp$ such that 
$|\vp| \leq M$, where $\|\cdot\|_\infty$ is the sup norm and
$M = \max\{|\psi| : \psi \in \supp x \cup \supp y_0\}$. \\

\noindent\underline{Claim}. Let $A$ be an admissible set such that
$\min\{|\vp| : \vp \in A\} \leq M$. 
If $\sum_{\vp\in A}|y_1(\vp)| = c$, and 
$\sum_{\vp\in A}|y_2(\vp)| = d$, then there exists an admissble set 
$B$ such that 
\begin{multline*} 
\min\{|\vp| : \vp \in A\} \leq \min\{|\vp| : \vp \in B\} \\
\leq \max\{|\vp| : \vp \in B\} \leq \max\{|\vp| : \vp \in A\}, 
\end{multline*}
$A \cap \supp y_0 \subseteq B$, 
and $\sum_{\vp \in B}|y_1(\vp)| \geq c + d$.\\

To prove the claim, let $N$ be such that $A \subseteq 
\cup_{|\vp| = N}T_\vp$, and $|A \cap T_\vp| \leq 1$ for all $\vp$ with
$|\vp| = N$. Then
$N \leq M$. Now, for each $\psi \in A \cap \supp y_2$, $\psi \in T_\vp$ 
for some $\vp$ with $|\vp| = N \leq M$. It follows that 
\[ \|y_1\chi_{T_\vp}\|_\infty = \|y_2\chi_{T_\vp}\|_\infty 
\geq |y_2(\psi)|.\]
Hence, there exists $\psi' \in T_\vp$ such that 
$|y_1(\psi')| \geq |y_2(\psi)|$. Now let
\[ B = (A \cap (\supp y_0 \cup \supp y_1)) 
\cup \{\psi' : \psi \in A \cap \supp y_2\}. \]
It is easy to see that the set $B$ satisfies the claim.\\
 
Suppose $\|x + x_1 + x_2\|_X/3 \geq  1 - \delta$. 
Let $x + x_1 + x_2 = x + 2y_0 + y_1 + y_2$ be normed by a 
sequence of admissible sets $A_1 \ll A_2 \ll \dots \ll A_k$. 
Denote by $\alpha = (a_i)\ki$, $\beta = (b_i)\ki$, 
$\gamma = (c_i)\ki$, and $\eta = (d_i)\ki$ 
respectively the 
sequences $(\sump|x(\vp)|)\ki$, $(\sump|y_0(\vp)|)\ki$, 
$(\sump|y_1(\vp)|)\ki$, and 
$(\sump|y_2(\vp)|)\ki$.  Now 
\[ \|\alpha + (\beta + \gamma) + (\beta + \eta)\|_2/3 
\geq \|x + x_1 + x_2\|_X/3 \geq  1 - \delta. \]
But $\|\alpha\|_2 \leq \|x\|_X \leq 1$. Similary, $\|\beta + \gamma\|_2$,
$\|\beta + \eta\|_2 \leq 1$. By Lemma \ref{trivial}, we obtain that
$\|\alpha - \beta - \gamma\|_2$, $\|\alpha - \beta - \eta\|_2$, and 
$\|\gamma - \eta\|_2$ are all $\leq \sqrt{18\delta}$. 
Let $j$ be the largest 
integer such that $a_j \neq 0$. Note that this implies 
$\supp x \cap A_j \neq \emptyset$; hence 
$(\supp y_1 \cup \supp y_2) \cap A_i = \emptyset$ for all $i < j$. 
Thus, $c_i = d_i = 0$ for all $i < j$. Now
\begin{equation} \label{b+d}
\|(b_{j+1}+d_{j+1},\dots,b_k+d_k)\|_2 \leq \|\alpha - \beta - \eta\|_2 \leq 
\sqrt{18\delta}.
\end{equation}
Moreover,
\begin{align}\label{beta}  
1 \geq \|x_2\|^2_X &= \|y_0 + y_2\|_X^2 \geq \|y_0\|_X^2 + \|y_2\|_X^2 \geq
\|\beta\|^2_2 + \|y_2\|_X^2 \notag\\
\implies \quad\quad \|\beta\|^2_2 &\leq 1 - \ep^2/9.
\end{align}
Hence
\begin{alignat*}{2}
3(1 - \delta) &\leq \|\alpha\|_2 + \|\beta + \gamma\|_2 + 
\|\beta + \eta\|_2 \leq 2 + \|\beta + \eta\|_2 & & \\
\implies \quad (1 - 3\delta)^2 &\leq \|\beta+\eta\|^2_2 & & \\
&= \|(b_1,\dots,b_{j-1},b_j+d_j)\|^2_2 &&\\
&\quad\quad + \|(b_{j+1}+d_{j+1},\dots,b_k+d_k)\|^2_2 & & \\
&\leq (\|(b_1,\dots,b_{j-1},b_j)\|_2 + d_j)^2 + 18\delta &   
&\text{by (\ref{b+d})}\\
&\leq (\|\beta\|_2 + d_j)^2 + 18\delta & & \\
&\leq ((1 - \ep^2/9)^{1/2} + d_j)^2 + 18\delta. & 
&\text{by (\ref{beta})}
\end{alignat*}
Therefore,
\begin{equation}\label{d}
d_j \geq (1 - 24\delta)^{1/2} - (1 - \ep^2/9)^{1/2}.
\end{equation}
Note that by the first part of the argument above, we also obtain that
\begin{equation}\label{beta+gamma}
\|\beta + \gamma\|_2 \geq 1 - 3 \delta.
\end{equation}
Since $A_j \cap \supp x \neq \emptyset$, we may apply the claim 
to obtain an admissible set $B$. Using the sequence of admissible sets
$A_1 \ll \dots \ll A_{j-1} \ll B \ll A_{j+1} \ll \dots \ll A_k$ to norm 
$x_1 = y_0 + y_1$ yields
\begin{align*}
1 &\geq \|y_0 + y_1\|_X^2 \geq 
\|(b_1,\dots,b_{j-1},b_j+c_j+d_j,b_{j+1}+c_{j+1},\dots,b_k+c_k)\|^2_2\\
&\geq 
\|(b_1,\dots,b_{j-1},b_j+c_j,b_{j+1}+c_{j+1},\dots,b_k+c_k)\|^2_2 + d^2_j
\\
&= \|\beta + \gamma\|^2_2 + d^2_j\\
&\geq (1 - 3\delta)^2 + [(1 - 24\delta)^{1/2} - (1 - \ep^2/9)^{1/2}]^2 
\quad \text{by (\ref{beta+gamma}) and (\ref{d})}.
\end{align*}
As the last expression is $> 1$ by (\ref{delta}), we have reached a 
contradiction.
\end{proof}

\noindent{\bf Remark.} The same method can be used to show 
that $X$ is $2$-$\beta$ with
the norm given in \cite{Sch}.\\

\begin{prop}
$(Y,\|\cdot\|_Y)$ is $2$-NUC.
\end{prop}

\begin{proof}
Let $(x_n)$ be an $\ep$-separated sequence in the unit ball of $Y$. Choose $\delta > 0$ so that 
\begin{align}
\label{Y1} \delta' = 12\delta + 2\sqrt{8\delta} &\leq \ep^2/18 \\
\intertext{and} 
\label{Y2} 1 - 2\delta - (2+\sqrt{8})\sqrt{\delta} &> \sqrt{1 - (\ep/3)^2}.
\end{align}
As in the proof of the previous proposition, it may be assumed that there exists a sequence $(y_n)^\infty_{n=0}$ in $Y$ such that $x_n = y_0 + y_n$, $\supp y_{n-1} \ll \supp y_n$ for all $n \in \N$, and $\|y_j\chi_{S_\vp}\|_\infty = \|y_k\chi_{S_\vp}\|_\infty$ whenever $|\vp| \leq M_i$ and $j, k > i$, where $M_i = \max\{|\psi| : \psi \in \supp y_i\}$ . We may also assume that $(\|y_n\|_Y)$ converges. Since $(y_n)^\infty_{n=1}$ is $\ep$-separated, $\eta = \lim \|y_n\|_Y \geq \ep/2$. The choice of $\delta'$ in (\ref{Y1}) guarantees that $4(\eta^2 - \delta')^{1/2} > 7\eta/2 \geq 3\eta + \sqrt{\delta'}$. Hence there exist $\eta_+ > \eta > \eta_- > \ep/3$ such that
\begin{equation}\label{Y3}
4\theta \geq 3\eta_+ + \sqrt{(\eta_+)^2 - (\eta_-)^2 + \delta'},
\end{equation}
where $\theta = \sqrt{(\eta_-)^2 - \delta'}$. We may now further assume that $\eta_+ \geq \|y_n\|_Y \geq \eta_-$ for all $n \in \N$. Now suppose that $\|x_m + x_n\|_Y/2 > 1 - \delta$ for all $m, n \in \N$. \\

\noindent\underline{Claim}. For all $m < n$ in $\N$, there exists an acceptable set $A$ such that $\sum_{\vp\in A}|y_i(\vp)| > \theta$ for $i = m, n $.\\

First observe that there are acceptable sets $A_1 \ll A_2 \ll \dots \ll A_k$ such that $\sum^k_{i=1}(\sum_{\vp\in A_i}|(2y_0+y_m+y_n)(\vp)|)^2 > 4(1-\delta)^2$. Let $\alpha = (a_i)^k_{i=1}$, $\beta = (b_i)^k_{i=1}$ and $\gamma = (c_i)^k_{i=1}$ be the sequences $(\sum_{\vp\in A_i}|y_j(\vp)|)^k_{i=1}$ for $j = 0, m, n$ respectively. Then $\|2\alpha+\beta+\gamma\|_2 > 2(1 - \delta)$ and $\|\alpha+\beta\|_2 \leq \|y_0+y_m\|_Y = \|x_m\|_Y \leq 1$. Similarly, $\|\alpha+\gamma\|_2 \leq 1$. It follows from the parallelogram law that $\|\beta - \gamma\|_2 < 4 - 4(1-\delta)^2 \leq 8 \delta$.  Note also that $\|\alpha + \beta\|_2 \geq \|2\alpha + \beta+ \gamma\|_2 - \|\alpha+\gamma\|_2 > 1 - 2\delta$. Similarly, $\|\alpha+\gamma\|_2 > 1 - 2\delta$. Let $j_1$, respectively, $j_2$, be the largest $j$ such that $a_j \neq 0$, respectively, $b_j \neq 0$. Since $\supp y_0 \cap A_{j_1} \neq \emptyset$, $b_1 = \dots = b_{j_1 - 1} = 0$. Similarly, $c_1 = \dots = c_{j_2 - 1} = 0$. Moreover, $j_1 \leq j_2$. Let us show that $j_1 < j_2$. For otherwise, $j_1 = j_2 = j$. Then
\begin{equation} \label{Y4}
|b_j - c_j| \leq \|\beta - \gamma\|_2 < \sqrt{8\delta}.
\end{equation}
Consider the set $A_j$. Choose $p \in \N \cup \{0\}$ such that $A_j \subseteq \cup_{|\vp|=p}S_\vp$ and $|A_j\cap S_\vp| \leq 1$ for all $\vp$ with $|\vp| = p$. Note that 
$p \leq M_0$. Let $G = \{\vp : |\vp| = p, A_j \cap S_\vp \cap \supp y_m \neq \emptyset\}$. If $\vp \in G$, $\|y_n\chi_{S_\vp}\|_\infty = \|y_m\chi_{S_\vp}\|_\infty$. Hence there exists $\psi_\vp \in S_\vp \cap \supp y_n$ such that $|y_n(\psi_\vp)| = \|y_m\chi_{S_\vp}\|_\infty$. It is easy to see that the set $B = \{\psi_\vp : \vp \in G\} \cup (A_j \cap \supp y_0)\ \cup (A_j \cap \supp y_n)$ is acceptable, and $\min\{|\vp| : \vp \in A_j\} \leq \min\{|\vp| : \vp \in B\}$. Hence $A_1 \ll \dots \ll A_{j-1} \ll B$. Thus
\begin{alignat*}{2}
1 &\geq \|x_n\|^2_Y = \|y_0 + y_n\|^2_Y 
\geq \sum^{j-1}_{i=1}|a_i|^2 + (\sum_{\vp\in B}|(y_0+y_n)(\vp)|)^2&&\\
&\geq \sum^{j-1}_{i=1}|a_i|^2 + (\sum_{\vp\in A_j}|y_0(\vp)| + \sum_{\vp\in A_j}|y_n(\vp)| + \sum_{\vp\in G}|y_n(\psi_\vp)|)^2 &&\\
&\geq \sum^{j-1}_{i=1}|a_i|^2 + (|a_j| + |c_j| + \sum_{\vp\in G}\|y_m\chi_{S_\vp}\|_\infty)^2 &&\\ 
&\geq \sum^{j-1}_{i=1}|a_i|^2 + (|a_j| + |c_j| + \sum_{\vp\in A_j}|y_m(\vp)|)^2 &&\\
&\geq \|(a_1,\dots,a_{j-1},a_j+b_j+c_j)\|^2_2 &&\\
&\geq \|(a_1,\dots,a_{j-1},a_j+b_j)\|^2_2 + |c_j|^2 &&\\
&\geq \|\alpha + \beta\|^2_2 + (|b_j| - \sqrt{8\delta})^2 &&\text{by (\ref{Y4})}\\
&> (1-2\delta)^2 + (\|\beta\|_2 - \sqrt{8\delta})^2.
\end{alignat*}
Therefore, $\|\beta\|_2 < (2+\sqrt{8})\sqrt{\delta}$. It follows that 
\begin{equation}\label{Y5}
\|\alpha\|_2 \geq \|\alpha+\beta\|_2 - \|\beta\|_2 > 1 - 2\delta - (2+\sqrt{8})\sqrt{\delta}.
\end{equation}
However, 
\begin{equation}\label{Y6}
 \|\alpha\|^2_2 \leq \|y_0\|^2_Y \leq \|x_m\|^2_Y - \|y_m\|^2_Y \leq 1 - (\eta_-)^2 < 1 - (\ep/3)^2. 
\end{equation}
Combining (\ref{Y5}) and (\ref{Y6}) with the choice of $\delta$ (\ref{Y2}) yields a contradiction. This shows that $j_1 < j_2$. Applying the facts that $\|\alpha + \beta\|_2 > 1 - 2\delta$ and $\|(b_{j_1},\dots,b_{j_2-1})\|_2 \leq \|\beta - \gamma\|_2 < \sqrt{8\delta}$, we obtain that 
\[ |b_{j_2}|^2 > (1 - 2\delta)^2 - (\|\alpha\|_2 + \sqrt{8\delta})^2 \geq (1 - 2\delta)^2 - (
\sqrt{1-(\eta_-)^2} + \sqrt{8\delta})^2 \geq \theta^2. \]
Similarly, 
\begin{align*}
(1 - 2\delta)^2 &< \|\alpha + \gamma\|^2_2 = \|\alpha\|^2_2 + |c_{j_2}|^2 + \|(c_{j_2+1},\dots,c_k)\|^2_2 \\
&\leq \|\alpha\|^2_2 + |c_{j_2}|^2 + \|\beta - \gamma\|^2_2 \\
&\leq 1 - (\eta_-)^2 + |c_{j_2}|^2 + 8\delta.
\end{align*}
Hence $|c_{j_2}| > \theta$. Thus the set $A = A_{j_2}$ satisfies the requirements of the claim.\\

Taking $m = 1$, $n = 2$, and $m = 2$, $n = 3$ respectively, we obtain acceptable sets $A$ and $A'$ from the claim. Since $A \cap \supp y_1 \neq \emptyset$, if $\vp \in A \cap \supp y_2$, $\vp \in S_{\vp'}$ for some $\vp'$ such that $|\vp'| \leq M_1$. This implies that there exists $\psi_\vp \in S_{\vp'}$ such that $|y_3(\psi_\vp)| = \|y_3\chi_{S_{\vp'}}\|_\infty = \|y_2\chi_{S_{\vp'}}\|_\infty \geq |y_2(\vp)|$. 
Let $q = \min\{|\vp| : \vp \in \supp y_3\}$ and $\Phi = \{\sigma \in T : |\sigma| = q\}$. For $\sigma \in \Phi$, define $s(\sigma) = |y_3(\psi_\vp)|$ if there exists $\vp \in A \cap \supp y_2$ such that $\psi_\vp \in S_\sigma$; otherwise, let $s(\sigma) = 0$. Also, let $t(\sigma) = |y_3(\vp)|$ if there exists $\vp \in A' \cap \supp y_3 \cap S_\sigma$; otherwise, let $t(\sigma) = 0$. Finally, let $r(\sigma) = \|y_3\chi_{S_\sigma}\|_\infty$ for all $\sigma \in \Phi$. Then $r(\sigma) \geq s(\sigma) \geq 0$ for all $\sigma \in \Phi$, $\sum_\sigma r(\sigma) \leq \|y_3\|_Y < \eta_+$, and $\sum_\sigma s(\sigma) > \theta$. Hence $\sum_\sigma(r(\sigma)-s(\sigma)) < \eta_+ - \theta$. Similarly, $\sum_\sigma(r(\sigma)-t(\sigma)) < \eta_+ - \theta$. Therefore, $\sum_\sigma|t(\sigma)-s(\sigma)| < 2(\eta_+ - \theta)$. Let $B$ be the set of all nodes in $A\cap \supp y_2$ that are comparable with some node in $A'\cap \supp y_3$. Then
\[ \sum_{\vp \in A \bs B}|y_2(\vp)| \leq \sum_{\vp \in A \bs B}|y_3(\psi_\vp)| \leq \sum_\sigma|t(\sigma)-s(\sigma)| < 2(\eta_+ - \theta).\]
Hence $\sum_{\vp \in B}|y_2(\vp)| > \theta - 2(\eta_+ - \theta) = 3\theta - 2\eta_+$. Now let $l = \min\{|\vp| : \vp \in A'\cap\supp y_2\}$. Divide $B$ into $B_1 = \{\vp \in B : |\vp| < l\}$ and $B_2 = \{\vp \in B : |\vp| \geq l\}$. Since $B_1$ and $A'\cap \supp y_2$ are acceptable sets such that $B_1 \ll A'\cap \supp y_2$, 
\begin{align*}
(\eta_+)^2 &> \|y_2\|^2_Y \geq (\sum_{\vp\in B_1}|y_2(\vp)|)^2 + (\sum_{\vp\in A'}|y_2(\vp)|)^2 \\ &> (\sum_{\vp\in B_1}|y_2(\vp)|)^2 + \theta^2 .
\end{align*}
Thus
\[
\sum_{\vp\in B_2}|y_2(\vp)| >  3\theta - 2\eta_+ - \sqrt{(\eta_+)^2 - \theta^2} .\]
Finally, since $B_2 \cup (A' \cap \supp y_2)$ is acceptable,
\begin{align*}
\eta_+ &> \|y_2\|_Y \geq \sum_{\vp\in B_2}|y_2(\vp)| + \sum_{\vp\in A' \cap \supp y_2}|y_2(\vp)| \\
&> 3\theta - 2\eta_+ - \sqrt{(\eta_+)^2 - \theta^2} + \theta.
\end{align*}
This contradicts inequality (\ref{Y3}).
\end{proof}

Before proceeding further, let us introduce some more notation. A branch in $T$ is a maximal subset of $T$ with respect to the partial order $\leq$.
If $\gamma$ is a branch in $T$, and $n \in \N \cup \{0\}$,
let $\vp^\gamma_n$ be the node of length $n$ in $\gamma$.
A collection of pairwise distinct branches is said to have 
separated at
level $L$ if for any pair of distinct branches $\gamma$ and $\gamma'$ 
in the collection, the nodes of length $L$ belonging to $\gamma$ and 
$\gamma'$ respectively are distinct. Finally, if $(\gamma_1,\dots,\gamma_k)$
is a sequence of pairwise distinct branches which have separated at a 
certain level $L$, say that a sequence of nodes $(\vp_1,\dots,\vp_k) \in
S(\gamma_1,\dots,\gamma_k;L)$ if $\vp_i \in T(\vp^{\gamma_i}_L)$, 
$1\leq i\leq k$. Let us note that
in this situation, $\|\chi_{\{\vp_i:1\leq i\leq k\}}\|_X = k$.

Suppose $|||\cdot|||$ is an equivalent norm on $X$ which is $2$-NUC. It may
be assumed that there exists $\ep > 0$ so that $\ep\|x\|_X \leq |||x||| \leq 
\|x\|_X$ for all $x \in X$. Let $\delta = \delta(2\ep) > 0$ be the number
obtained from the definition of $2$-NUC for the norm $|||\cdot|||$.

\begin{prop}\label{main}
Let $n \in \N\ \cup \{0\}$. Then there are pairwise incomparable nodes 
$\vp_1,\dots,\vp_{2^n}$ such that whenever $\gamma_i, \gamma'_i$ 
are distinct branches passing through $\vp_i$, $1 \leq i \leq 2^n$, and 
$\{\gamma_i,\gamma'_i: 1\leq i \leq 2^n\}$ have separated at level $L$, 
there is a sequence of nodes $(\psi_1,\dots,\psi_{2^{n+1}}) \in 
S(\gamma_1,\gamma'_1,\dots,\gamma_{2^n},\gamma'_{2^n};L)$ satisfying
$|||\chi_{\{\psi_i:1\leq i \leq 2^{n+1}\}}||| \leq (2(1-\delta))^{n+1}$.
\end{prop}

\begin{proof}
Assume that $n$ is the first non-negative integer where the proposition 
fails. Let $\vp_1,\dots,\vp_{2^{n-1}}$ be the nodes obtained by applying
the proposition for the case $n-1$. (If $n =0$, begin the argument with 
any node $\vp_1$.) For each $i$, $1 \leq i \leq 2^{n-1}$, let $\psi_{2i-1,1}$
and $\psi_{2i,1}$ be a pair of incomparable nodes in $T_{\vp_i}$.
(If $n = 0$, let $\psi_{1,1}$ be any node in $T_{\vp_1}$.)
Since the proposition fails for the nodes $\psi_{1,1},\dots,\psi_{2^n,1}$,
there are distinct branches $\gamma_{i,1}$, $\gamma'_{i,1}$ passing through
$\psi_{i,1}$, $1 \leq i \leq 2^n$, and a number $L_1$ so that 
$\{\gamma_{i,1},\gamma'_{i,1}: 1\leq i \leq 2^n\}$ have separated at level 
$L_1$, but $|||\chi_{\{\xi_i:1\leq i\leq 2^{n+1}\}}||| > 
(2(1-\delta))^{n+1}$ for
any sequence of nodes $(\xi_1,\dots,\xi_{2^{n+1}}) \in 
S(\gamma_{1,1},\gamma'_{1,1},\dots,\gamma_{2^n,1},\gamma'_{2^n,1};L_1)$.
However, since the proposition holds for the nodes 
$\vp_1,\dots,\vp_{2^{n-1}}$, we obtain a sequence of nodes 
$(\xi_{1,1},\dots,\xi_{2^n,1}) \in S(\gamma'_{1,1},\dots,\gamma'_{2^n,1};L_1)$
such that 
\[ |||\chi_{\{\xi_{i,1}:1\leq i\leq 2^n\}}||| \leq (2(1-\delta))^n. \]
(Note that the preceding statement holds trivially if $n = 0$.)
For each $i$, choose a node $\psi_{i,2}$ in $\gamma_{i,1}$ such that
$|\psi_{i,2}| > L_1$. Then $\psi_{2i-1,2}$ 
and $\psi_{2i,2}$ is a pair of incomparable nodes in $T_{\vp_i}$, and the 
argument may be repeated. (If $n = 0$, repeat the argument using the 
node $\psi_{1,2}$.) Inductively, we thus obtain sequences of 
branches 
$(\gamma_{1,r},\gamma'_{1,r},\dots,\gamma_{2^n,r},\gamma'_{2^n,r})
^\infty_{r=1}$, 
a sequence of numbers $L_1 < L_2 < \dots$, and sequences of nodes 
$(\xi_{1,r},\dots,\xi_{2^n,r})^\infty_{r=1}$ such that
\begin{enumerate}
\item the branches $\{\gamma_{i,r},\gamma'_{i,r}:1 \leq i \leq 2^{n}\}$ have 
separated at level $L_r$, $r \geq 1$,
\item \label{lowerbound}
$|||\chi_{\{\xi_i:1 \leq i \leq 2^{n+1}\}}||| > (2(1-\delta))^{n+1}$
for any sequence of nodes 
\[ (\xi_1,\dots,\xi_{2^{n+1}}) \in 
S(\gamma_{1,r},\gamma'_{1,r},\dots,\gamma_{2^n,r},\gamma'_{2^n,r};L_r), \]
\item \label{xr} $(\xi_{1,r},\dots,\xi_{2^n,r}) \in 
S(\gamma'_{1,r},\dots,\gamma'_{2^n,r};L_r)$, and 
\[ |||\chi_{\{\xi_{i,r}:1 \leq i \leq 2^n\}}||| 
\leq (2(1-\delta))^n, \quad r \geq 1, \]
\item $\xi_{i,r} \in T(\vp^{\gamma_{i,s}}_{L_s})$ whenever $r > s$, and  
$1\leq i \leq 2^n$.
\end{enumerate}
It follows that if $r > s$, then
\begin{equation}\label{rs}
(\xi_{1,r},\xi_{1,s},\dots,\xi_{2^n,r},\xi_{2^n,s}) \in 
S(\gamma_{1,s},\gamma'_{1,s},\dots,\gamma_{2^n,s},\gamma'_{2^n,s};L_s).
\end{equation}
Let $x_r = 
(2(1-\delta))^{-n}\chi_{\{\xi_{i,r}:1\leq i\leq 2^n\}}$, $r \geq 1$. By item
\ref{xr}, $|||x_r||| \leq 1$. Moreover, because of (\ref{rs}), if $r > s$,
then
\[ |||x_r -  x_s||| \geq \ep\|x_r - x_s\|_X = 2^{n+1}\ep/(2(1-\delta))^n 
\geq 2\ep .\]
Thus $(x_r)$ is $2\ep$-separated in the norm $|||\cdot|||$. By the choice of
$\delta$, there are $r > s$ such that $|||x_r + x_s|||/2 \leq 1 - \delta$.
Therefore, 
$|||\chi_{\{\xi_{1,r},\xi_{1,s},\dots,\xi_{2^n,r},\xi_{2^n,s}\}}||| 
\leq (2(1-\delta))^{n+1}$. But this contradicts item \ref{lowerbound} and 
condition (\ref{rs}).
\end{proof}

\begin{thm} \label{th}
There is no equivalent $2$-NUC norm on $X$. 
\end{thm}

\begin{proof}
In the notation of the statement of Proposition \ref{main}, we obtain for
each $n$, nodes $\psi_1,\dots,\psi_{2^{n+1}}$ such that 
$|||\chi_{\{\psi_i:1\leq i \leq 2^{n+1}\}}||| \leq (2(1-\delta))^{n+1}$, and 
$\|\chi_{\{\psi_i:1\leq i \leq 2^{n+1}\}}\|_X = 2^{n+1}$. Hence $|||\cdot|||$ 
cannot be an equivalent norm on $X$.
\end{proof}

The proof that the space $Y$ has no equivalent $1$-$\beta$ norm follows along similar lines. Suppose that $|||\cdot|||$ is an equivalent $1$-$\beta$ norm on $Y$. We may assume that $\ep\|\cdot\|_Y \leq |||\cdot||| \leq \|\cdot\|_Y$ for some $\ep > 0$. Let $\delta = \delta(\ep)$ be the constant obtained from the definition of $1$-$\beta$ for the norm $|||\cdot|||$. Let $n \in \N \cup \{0\}$ and denote the set $\{\vp \in T : |\vp| = n\}$ by $\Phi$.

\begin{prop}
For any $m$, $0 \leq m \leq n$, any subset $\Phi'$ of $\Phi$ with $|\Phi'| = 2^m$, and any $p \in \N$, there exists an acceptable set $A \subseteq \cup_{\vp\in\Phi'}S_\vp$ such that $|A| = 2^m$, $\min\{|\vp| : \vp \in A\} \geq p$, and $|||\chi_A||| \leq 2^m(1-\delta)^m$.
\end{prop}

\begin{proof}
The case $m = 0$ is trivial. Suppose the proposition holds for some $m$, $0 \leq m < n$. Let $\Phi' \subseteq \Phi$, $|\Phi'| = 2^{m+1}$, and let $p \in \N$. Divide $\Phi'$ into disjoint subsets $\Phi_1$ and $\Phi_2$ such that $|\Phi_1| = |\Phi_2| = 2^m$. By the inductive hypothesis, there exist acceptable sets $B$ and $C_j$, $j \in \N$, such that $B \subseteq \cup_{\vp\in \Phi_1}S_\vp$, $|B| = 2^m$, $\min\{|\vp| : \vp \in B\} \geq p$, and $|||\chi_B||| \leq 2^m(1-\delta)^m$; and also $C_j \subseteq \cup_{\vp\in \Phi_2}S_\vp$, $|C_j| = 2^m$, $\min\{|\vp| : \vp \in C_1\} \geq p$, $C_{j} \ll C_{j+1}$, and $|||\chi_{C_j}||| \leq 2^m(1-\delta)^m$ for all $j \in \N$. It is easily verified that the sequence $(2^{-m}(1-\delta)^{-m}\chi_{C_j})$ is $\ep$-separated and has norm bounded by $1$ with respect to $|||\cdot|||$. It follows that there exists $j_0$ such that $2^{-m}(1-\delta)^{-m}|||\chi_B+\chi_{C_{j_0}}||| \leq 2(1-\delta)$. The induction is completed by taking $A$ to be $B \cup C_{j_0}$.
\end{proof}

Using the same argument as in Theorem \ref{th}, we obtain

\begin{thm}
There is no equivalent $1$-$\beta$ norm on $Y$.
\end{thm}

We close with the obvious problem.\\

\noindent{\bf Problem.} For $k \geq 3$, can every $k$-NUC Banach space, respectively $k$-$\beta$ Banach space, be equivalently renormed to be $(k-1)$-$\beta$, respectively, $k$-NUC?



\begin{thebibliography}{99}
\bibitem{Sch} {\sc Walter Schachermayer}, 
{\em The class of Banach spaces which do not
have $c_0$ as a spreading model is not $L^2$-hereditary}, 
Ann.\ Inst.\ H.\ Poincar\'{e} {\bf 19}(1983), 1-8.
\end{thebibliography}
\end{document}